\documentclass[11 pt]{amsart}

\usepackage{amsfonts}
\usepackage{amssymb}
\usepackage{amsmath}
\usepackage{latexsym}
\usepackage{mathrsfs}
\usepackage{amsthm}
\usepackage{verbatim}

\newtheorem{thrm}{Theorem}[section]
\newtheorem{lem}[thrm]{Lemma}
\newtheorem{cor}[thrm]{Corollary}
\newtheorem{prop}[thrm]{Proposition}

\theoremstyle{definition}
\newtheorem{defn}[thrm]{Definition}

\newtheorem{exmple}[thrm]{Example}

\newtheorem{ques}[thrm]{Question}

\begin{document}

\newcommand{\vol}{\mathrm{vol}}

\newcommand{\Supp}{\mathrm{Supp}}

\newcommand{\Sing}{\mathrm{Sing}}

\newcommand{\ord}{\mathrm{ord}}

\newcommand{\Pic}{\mathrm{Pic}}

\newcommand{\codim}{\mathrm{codim}}

\newcommand{\mult}{\mathrm{mult}}

\title{Numerical Triviality and Pullbacks}
\author{Brian Lehmann}
\thanks{The author is supported by NSF Award 1004363.}
\address{Department of Mathematics, Rice University \\
Houston, TX 77005}
\email{blehmann@rice.edu}

\begin{abstract}
Let $f: X \to Z$ be a surjective morphism of smooth complex projective varieties with connected fibers.  Suppose that $L$ is a pseudo-effective divisor that is $f$-numerically trivial.  We show that there is a divisor $D$ on $Z$ such that $L \equiv f^{*}D$.
\end{abstract}

\maketitle

\section{Introduction}

We consider the following question:

\begin{ques} \label{mainquestion}
Let $f: X \to Z$ be a surjective morphism of smooth complex projective varieties with connected fibers.  Suppose that $L$ is a pseudo-effective $\mathbb{R}$-Cartier divisor that is numerically trivial on the fibers of $f$.  Is $L$ numerically equivalent to the pull-back of a divisor on $Z$?
\end{ques}

When $L$ is not pseudo-effective the answer is an emphatic ``no.''  Thus it is perhaps surprising that there is a positive answer for pseudo-effective divisors.  The analogue of Question \ref{mainquestion} for $\mathbb{Q}$-linear equivalence is well understood, with the most general statements due to \cite{nakayama04}.  Our goal is to show that similar theorems hold true in the numerical case.

The most restrictive situation is to ask that $L$ be numerically trivial on \emph{every} fiber of $f$.  In this case $L$ is actually numerically equivalent to the pullback of a divisor on $Z$:

\begin{thrm} \label{numtrivialeveryfiber}
Let $f: X \to Z$ be a surjective morphism with connected fibers from a normal complex projective variety $X$ to a $\mathbb{Q}$-factorial normal complex projective variety $Y$.  Suppose that $L$ is a pseudo-effective $\mathbb{R}$-Cartier divisor that is $f$-numerically trivial.  Then there is some $\mathbb{R}$-Cartier divisor $D$ on $Z$ such that $L \equiv f^{*}D$.
\end{thrm}

Again, the pseudo-effectiveness of $L$ is crucial: the dimension of the space of divisors that are $f$-numerically trivial will generally be larger than the dimension of $N^{1}(Z)$.

For applications it is more useful to require that $L$ be numerically trivial only on a \emph{general} fiber of $f$.  To handle this case we need a systematic way of discounting the non-trivial behavior along special fibers.  For surfaces the behavior of special fibers is captured by the Zariski decomposition.  The analogous construction in higher dimensions is the divisorial Zarkiski decomposition of \cite{nakayama04} and \cite{boucksom04}.  Given a pseudo-effective $\mathbb{R}$-Cartier divisor $L$, the divisorial Zariski decomposition
\begin{equation*}
L = P_{\sigma}(L) + N_{\sigma}(L)
\end{equation*}
expresses $L$ as the sum of a ``movable part'' $P_{\sigma}(L)$  and a ``fixed part''  $N_{\sigma}(L)$ (see Definition \ref{divzardecdef}).  The following theorem is a numerical analogue of \cite{nakayama04} V.2.26 Corollary.  

\begin{thrm}  \label{numtrivialgeneralfiber}
Let $f: X \to Z$ be a surjective morphism of normal complex projective varieties with connected fibers.  Suppose that $L$ is a pseudo-effective $\mathbb{R}$-Cartier divisor such that $L|_{F} \equiv 0$ for a general fiber $F$ of $f$.  Then there is a smooth birational model $\phi: Y \to X$, a map $g: Y \to Z'$ birationally equivalent to $f$, and an $\mathbb{R}$-Cartier divisor $D$ on $Z'$ such that $P_{\sigma}(\phi^{*}L) \equiv P_{\sigma}(g^{*}D)$.
\end{thrm}

The most general situation is to ask that $\nu(L|_{F}) = 0$ on a general fiber $F$, where $\nu$ denotes the numerical dimension of \cite{nakayama04} and \cite{bdpp04}.  Theorem \ref{numtrivialgeneralfiber} also holds with this weaker condition on $L$.

To apply these results, one must find a morphism $f: X \to Z$ such that $L$ is numerically trivial along the fibers.  \cite{8authors}, \cite{eckl05}, and \cite{lehmann11} show that such maps can be constructed by taking the quotient of $X$ by subvarieties along which $L$ is numerically trivial.  In fact, there is a maximal fibration such that $L$ is numerically trivial (properly interpreted) along the fibers.  Theorems \ref{numtrivialeveryfiber} and \ref{numtrivialgeneralfiber} pair naturally with the reduction map theory developed in these papers.

\subsection{Outline}
The first step in the proofs is to construct a candidate divisor $D$ on the base $Z$.  We accomplish this by cutting down by very ample divisors to obtain a generically finite map $f: W \to Z$.  We then use facts about finite maps to ``push down'' $L|_{W}$ to obtain $D$.

The next step is to compare $L$ and $f^{*}D$ using numerical analogues of results of \cite{nakayama04}.  This is achieved by cutting down to the surface case.  An important conceptual point is that a numerical class in $N^{1}(X)$ is determined ``in codimension 1'': a divisor class is determined by its intersections with curves avoiding a countable union of closed subsets of codimension at least $2$.  

Section \ref{prelimsection} is devoted to preliminaries: surfaces, the divisorial Zariski decomposition, the numerical dimension, and exceptional divisors.
Section \ref{genfinsection} analyzes generically finite maps.
Sections \ref{everyfibersection} and \ref{generalfibersection} prove Theorems \ref{numtrivialeveryfiber} and \ref{numtrivialgeneralfiber} respectively.

I thank B.~Bhatt, A.M.~Fulger, and Y.~Gongyo for some helpful conversations.

\section{Preliminaries} \label{prelimsection}

We work over the base field $\mathbb{C}$.  All varieties are irreducible and reduced.

\subsection{Notation}

We will use the standard notations $\sim, \sim_{\mathbb{Q}}, \sim_{\mathbb{R}}$, and $\equiv$ to denote respectively linear equivalence, $\mathbb{Q}$-linear equivalence, $\mathbb{R}$-linear equivalence, and numerical equivalence.

\begin{defn}
Suppose that $f: X \to Z$ is a morphism of normal projective varieties.  We say that
\begin{itemize}
\item A curve $C$ on $X$ is $f$-vertical if $f^{*}H \cdot C = 0$ for some ample Cartier divisor $H$ on $Z$ (and similarly for a curve class $\alpha \in N_{1}(X)$).
\item An $\mathbb{R}$-Cartier divisor $L$ on $X$ is $f$-numerically trivial if it has vanishing intersection with every $f$-vertical curve class.
\end{itemize}
\end{defn}

Suppose that $f: X \to Z$ is a surjective morphism of normal projective varieties.  By \cite{raynaud72} there is a birational model $\psi: T \to Z$ such that for the normalization $W$ of the main component of $X \times_{Z} T$ the induced map $g: W \to T$ is flat.  We say that $g: W \to T$ is a flattening of $f$.

\subsection{Surfaces} \label{surfacesection}

We begin by considering Question \ref{mainquestion} for surfaces.  For surfaces, the Zariski decomposition is the key tool.

\begin{thrm}[\cite{zariski64},\cite{fujita79}]
Let $S$ be a smooth projective surface and let $L$ be a pseudo-effective $\mathbb{R}$-Cartier divisor on $S$.  There is a unique decomposition
\begin{equation*}
L = P +N
\end{equation*}
where $P$ is a nef divisor and $N$ is an effective divisor satisfying
\begin{enumerate}
\item $P \cdot N = 0$.
\item If $N \neq 0$, the intersection matrix defined by the components of $N$ is negative definite.
\end{enumerate}
\end{thrm}

We will use this theorem to describe how a divisor $L$ behaves along special fibers of an $L$-trivial morphism.  The following lemma is well-known.

\begin{lem} \label{surfacelemma}
Let $f: S \to C$ be a surjective morphism with connected fibers from a smooth projective surface $S$ to a smooth projective curve $C$.  Let $L$ be a pseudo-effective $\mathbb{R}$-Cartier divisor on $S$ such that $L \cdot F = 0$ for a general fiber $F$ of $f$.
\begin{enumerate}
\item If $L \cdot D = 0$ for every $f$-vertical curve $D$, then $L$ is nef.
\item If $L \cdot D \neq 0$ for some $f$-vertical curve $D$ contained in a fiber $F_{0}$, then there is an $f$-vertical curve $G$ contained in $F_{0}$ satisfying $L \cdot G < 0$.
\end{enumerate}
\end{lem}

We recall the proof for convenience.

\begin{proof}
Let $L = P + N$ be the Zariski decomposition of $L$.  Since $P$ is nef and $P \cdot F \leq L \cdot F = 0$, $P$ has vanishing intersection with every $f$-vertical curve.  Note that $N$ is an (effective) $f$-vertical curve since $N \cdot F = L \cdot F = 0$.

We first show (2).  By assumption $N$ must have some components contained in $F_{0}$.  Recall that the self-intersection matrix of the components of $N$ is negative-definite.  In fact, since $f$-vertical curves in different fibers do not intersect, the same is true just for the components contained in $F_{0}$.  Thus, there is an effective curve $G$ supported on $\Supp(N) \cap \Supp(F_{0})$ with $0 > N \cdot G = L \cdot G$.  The same argument shows that in (1) we must have $N=0$ so that $L$ is nef.
\end{proof}

The following is a special case of a theorem of \cite{peternell08}.  

\begin{lem}[cf.~\cite{peternell08},Theorem 6.8] \label{relativebdpp}
Let $f: S \to C$ be a surjective morphism with connected fibers from a smooth projective surface $S$ to a smooth projective curve $C$.  Suppose that $L$ is an $f$-numerically trivial nef $\mathbb{R}$-Cartier divisor.  Then $L \equiv \alpha F$ for some $\alpha \geq 0$ where $F$ denotes a general fiber of $f$.
\end{lem}

\begin{proof}
Suppose the theorem fails.  There is a divisor $D$  such that $D \cdot L < 0$ and $D \cdot F > 0$.  The latter condition implies that $D$ is $f$-big so that $D + f^{*}H$ is pseudo-effective for some ample divisor $H$ on $C$.  But $(D + mf^{*}H) \cdot L < 0$, a contradiction.
\end{proof}

\begin{cor} \label{numtrivsurfaces}
Let $f: S \to C$ be a surjective morphism from an irreducible projective surface $S$ to a smooth projective curve $C$ with connected fibers.  Suppose that $L$ is a pseudo-effective $\mathbb{R}$-Cartier divisor on $S$ such that $L \cdot C = 0$ for every $f$-vertical curve $C$.  Then $L \equiv f^{*}D$ for some $\mathbb{R}$-Cartier divisor $D$ on $T$.
\end{cor}

\begin{proof}
When $S$ is smooth this follows from Lemmas \ref{surfacelemma} and \ref{relativebdpp}.  In general we pass to a resolution $\phi: S' \to S$.  Applying the smooth case to $\phi^{*}L$ we find a divisor $D$ such that $\phi^{*}L \equiv (f \circ \phi)^{*}D$.  Thus $L \equiv f^{*}D$.
\end{proof}

\subsection{Divisorial Zariski decompositions}

We next recall the divisorial Zariski decomposition.  This notion was introduced by \cite{nakayama04} and \cite{boucksom04} as a higher-dimensional analogue of the Zariski decomposition for surfaces.

\begin{defn}
Let $X$ be a smooth projective variety and let $L$ be a pseudo-effective $\mathbb{R}$-Cartier divisor on $X$.  Fix an ample divisor $A$ on $X$.  Given a prime divisor $\Gamma$ on $X$, we define
\begin{equation*}
\sigma_{\Gamma}(L) = \lim_{\epsilon \to 0} \min \{\mult_{\Gamma}(L') | L' \geq 0 \textrm{ and } L' \sim_{\mathbb{Q}} L+\epsilon A\}.
\end{equation*}
This definition is independent of the choice of $A$.
\end{defn}

\cite{nakayama04} shows that for any pseudo-effective divisor $L$ there are only finitely many prime divisors $\Gamma$ with $\sigma_{\Gamma}(L) > 0$.  Thus we can make the following definition.

\begin{defn} \label{divzardecdef}
Let $X$ be a smooth projective variety and let $L$ be a pseudo-effective $\mathbb{R}$-Cartier divisor.
We define:
\begin{equation*}
N_{\sigma}(L) = \sum \sigma_{\Gamma}(L) \Gamma \qquad \qquad P_{\sigma}(L) = L - N_{\sigma}(L)
\end{equation*}
The decomposition $L = P_{\sigma}(L) + N_{\sigma}(L)$ is called the divisorial Zariski
decomposition of $L$.
\end{defn}

We need the following properties of the divisorial Zariski decomposition.

\begin{lem}[\cite{nakayama04}, III.1.4 Lemma, V.1.3 Theorem, and III.2.5 Lemma] \label{divzarprop}
Let $X$ be a smooth projective variety and let $L$ be a pseudo-effective $\mathbb{R}$-Cartier divisor.  Then
\begin{enumerate}
\item $N_{\sigma}(L)$ is effective.
\item For any prime divisor $\Gamma$ of $X$, the restriction $P_{\sigma}(L)|_{\Gamma}$ is pseudo-effective.
\item If $\phi: Y \to X$ is a birational map from a smooth projective variety $Y$, then $N_{\sigma}(\phi^{*}L) \geq \phi^{*}N_{\sigma}(L)$.
\end{enumerate}
\end{lem}

The following is a numerical analogue of \cite{nakayama04} III.5.2 Lemma.

\begin{lem}[\cite{lehmann11}, Lemma 4.4] \label{relnotpsef}
Let $f: X \to Z$ be a surjective morphism from a smooth projective variety to a normal projective variety with connected fibers.  Suppose that $L$ is an $\mathbb{R}$-Cartier divisor such that $L|_{F}  \equiv 0$ on the general fiber $F$ of $f$.  If $\Theta$ is a prime divisor on $Z$ such that $L|_{F} \not \equiv 0$ for a general fiber $F$ over $\Theta$, then there is some prime divisor $\Gamma$ on $X$ such that $f(\Gamma)=\Theta$ and $L|_{\Gamma}$ is not pseudo-effective.
\end{lem}

\begin{proof}
The surface case is Lemma \ref{surfacelemma} (2).  The general case is proved by cutting down by general very ample divisors on $X$ and $Z$ to reduce to the surface case.
\end{proof}

\begin{cor} \label{psigmarelativelytrivial}
Let $f: X \to Z$ be a surjective morphism from a smooth projective variety to a normal projective variety with connected fibers.  Suppose that $L$ is an $\mathbb{R}$-Cartier divisor such that $L|_{F}  \equiv 0$ on the general fiber $F$ of $f$.  Then there is a subset $V \subset Z$ that is a countable union of closed sets of codimension $2$ such that $P_{\sigma}(L)|_{F} \equiv 0$ for every fiber $F$ not lying above $V$.
\end{cor}

\begin{proof}
Since $L \geq P_{\sigma}(L)$, we see that $P_{\sigma}(L)|_{F} \equiv 0$ for a general fiber $F$ of $f$.  The conclusion follows from Lemma \ref{relnotpsef} combined with Lemma \ref{divzarprop} (2).
\end{proof}

\subsection{Numerical dimension}

The numerical dimension is a numerical measure of positivity that is closely related to the divisorial Zariski decomposition.  We will use the definition of \cite{nakayama04}.

\begin{defn}
Let $X$ be a normal projective variety and let $L$ be a pseudo-effective $\mathbb{R}$-Cartier divisor on $X$.  For an ample divisor $A$ set
\begin{equation*}
\nu(L,A) = \max \left\{ k \in \mathbb{Z}_{\geq 0} \left| \limsup_{m \to \infty} \frac{h^{0}(X,\lfloor mL \rfloor + A)}{m^{k}} > 0 \right. \right\}
\end{equation*}
and define
\begin{equation*}
\nu(L) = \max_{A \textrm{ ample}} \{ \nu(L,A) \}.
\end{equation*}
\end{defn}

\begin{lem}
Let $X$ be a normal projective variety and let $L$ be a pseudo-effective $\mathbb{R}$-Cartier divisor.  Then
\begin{enumerate}
\item If $\phi: Y \to X$ is a birational map from a normal projective variety $Y$, then $\nu(\phi^{*}L) = \nu(L)$.
\item If $X$ is smooth then $\nu(L) = 0$ iff $P_{\sigma}(L) \equiv 0$.
\end{enumerate}
\end{lem}

\subsection{Exceptional divisors}

\begin{defn}
Let $f: X \to Z$ be a surjective morphism of normal projective varieties.  An $\mathbb{R}$-Cartier divisor $E$ on $X$ is
\begin{itemize}
\item $f$-vertical if no component of $\Supp(E)$ dominates $Z$.
\item $f$-horizontal otherwise.
\end{itemize}
\end{defn}

We next identify two different ways a divisor can be ``exceptional'' for a morphism.

\begin{defn}
Let $f: X \to Z$ be a surjective morphism of normal projective varieties.  An $f$-vertical $\mathbb{R}$-Cartier divisor $E$ on $X$ is
\begin{itemize}
\item $f$-exceptional if every component $E_{i}$ of $\Supp(E)$ satisfies $$\codim f(E_{i}) \geq 2.$$
\item $f$-degenerate if for every prime divisor $\Theta \subset Z$ there is a prime divisor $\Gamma \subset X$ with $f(\Gamma)=\Theta$ and $\Gamma \not \subset \Supp(E)$.
\end{itemize}
\end{defn}

\begin{lem}[cf.~\cite{nakayama04} III.5.8 Lemma] \label{makefdegen}
Let $f: X \to Z$ be a surjective morphism of smooth projective varieties
with connected fibers.  Suppose that $L$ is an effective $f$-vertical $\mathbb{R}$-Cartier divisor.  There is an effective $\mathbb{R}$-Cartier divisor $D$ on $Z$ and an effective $f$-exceptional divisor $E$ on $X$ such that
\begin{equation*}
L + E = f^{*}D + F
\end{equation*}
where $F$ is an effective $f$-degenerate divisor.
\end{lem}

\begin{proof}
Write $L = L_{exc} + L_{ft}$ where $L_{exc}$ is the part of $L$ supported on the $f$-exceptional locus and $L_{ft}$ is the rest.

Suppose that $\{ \Theta_{i} \}_{i=1}^{r}$ are the prime divisors on $Z$ contained in the image $f(\Supp(L))$.  To $\Theta_{i}$ assign the constant $\beta_{i}$ defined by
\begin{equation*}
\beta_{i} = \min_{\substack{ \Gamma \textrm{ a prime divisor on }X \\ \textrm{ with } f(\Gamma) = \Theta_{i}} } \frac{\mult_{\Gamma}(L_{ft})}{\mult_{\Gamma}(f^{*}\Theta_{i})}.
\end{equation*}
Set $D = \sum_{i=1}^{n} \beta_{i} \Theta_{i}$.  There is some $f$-exceptional divisor $E$ such that $F := L + E - f^{*}D$ is an effective $f$-vertical divisor.  Furthermore, by construction we know that for each $\Theta$ there is a prime component $\Gamma$ of $X$ such that $f(\Gamma) = \Theta$ but $\Gamma \not \subset \Supp(F)$.
\end{proof}

As demonstrated by Nakayama, the divisorial Zariski decomposition gives a useful language for understanding $f$-degenerate divisors.  

\begin{lem}[\cite{gl11}, Lemma 2.16] \label{fdegen}
Let $f: X \to Z$ be a surjective morphism from a smooth projective variety to a normal projective variety and let $D$ be an effective $f$-degenerate $\mathbb{R}$-Cartier divisor.  For any pseudo-effective $\mathbb{R}$-Cartier divisor $L$ on $Z$ we have $D \leq N_{\sigma}(f^{*}L + D)$.
\end{lem}

\section{Generically Finite Maps} \label{genfinsection}

In this section we study the behavior of divisors over generically finite morphisms.  Such morphisms are a composition of a birational map and a finite map and can be understood by addressing each separately.  The following lemma is a well-known consequence of the Negativity of Contraction lemma.

\begin{lem} \label{birationallem}
Let $f: X \to Z$ be a birational morphism from a normal projective variety $X$ to a $\mathbb{Q}$-factorial normal projective variety $Z$.  Suppose that $L$ is an $\mathbb{R}$-Cartier divisor on $X$ such that $L$ is $f$-numerically trivial.  Then $L \equiv f^{*}D$ for some $\mathbb{R}$-Cartier divisor $D$ on $Z$.
\end{lem}

\begin{proof}
Suppose that $\phi: X' \to X$ is a resolution of $X$.  Note that $\phi^{*}L$ is $(f \circ \phi)$-numerically trivial.  If we can find a $D$ so that $\phi^{*}L \equiv (f \circ \phi)^{*}D$ then also $L \equiv f^{*}D$.  Thus we may assume that $X$ is smooth.

Write $L \equiv f^{*}f_{*}L + E$ for some $f$-exceptional divisor $E$.  We can write $E = E^{+} - E^{-}$ for some effective divisors $E^{+},E^{-}$ with disjoint support.  Since $L$ is $f$-numerically trivial, we have $E^{+} \equiv_{f} E^{-}$.  However, \cite{nakayama04} III.5.1 Lemma implies that any non-zero $f$-exceptional divisor has negative intersection with some $f$-vertical curve that deforms to cover a divisor.  Since the supports of $E^{+}$ and $E^{-}$ do not share any component we must have $E^{+} = E^{-} = 0$.
\end{proof}

\begin{lem} \label{finitepullback}
Let $f: X \to Z$ be a surjective finite morphism of normal projective varieties and let $L$ be an $\mathbb{R}$-Cartier divisor on $X$.  Let $\{ T_{i} \}_{i = 1}^{k}$ be a collection of irreducible curves on $X$.  Suppose that there are constants $\alpha_{i}$ such that
\begin{equation*}
L \cdot C = (\deg f|_{C}) \alpha_{i}
\end{equation*}
for every curve $C$ on $X$ with $f(C)=f(T_{i})$.  Then there is an $\mathbb{R}$-Cartier divisor $D$ on $Z$ such that $L \cdot T_{i} = f^{*}D \cdot T_{i}$ for every $i$.  In particular, if the numerical classes of the $T_{i}$ span $N_{1}(X)$ then $L \equiv f^{*}D$. 
\end{lem}

\begin{proof}
Let $h: W \to Z$ denote the Galois closure of $f$ with Galois group $G$.  We let $p: W \to X$ denote the map to $X$.  We first show that the $\mathbb{R}$-Cartier divisor
\begin{equation*}
L_{G} := \frac{1}{|G|}\sum_{g \in G} g(p^{*}L)
\end{equation*}
is numerically equivalent to $h^{*}D$ for some $\mathbb{R}$-Cartier divisor $D$ on $Z$.  For a positive integer $m$, let $L_{m} = \sum_{g \in G} g(\lfloor mp^{*}L \rfloor)$.  Since $L_{m}$ is $G$-invariant, we can find a Cartier divisor $D_{m}$ on $Z$ such that $h^{*}D_{m} = L_{m}$.  Note that the numerical classes of $\frac{1}{m|G|}D_{m}$ converge; choose $D$ to be an $\mathbb{R}$-Cartier divisor representing this class.  Then
\begin{equation*}
h^{*}D \equiv \lim_{m \to \infty} \frac{1}{m|G|} D_{m} = L_{G}.
\end{equation*}
Note that if $C$ is a curve on $W$ such that $p(C) = T_{i}$ then $L_{G} \cdot C = p^{*}L \cdot C$ by the assumption on the intersection numbers of $L$.  Thus $L \cdot T_{i} = f^{*}D \cdot T_{i}$ for each $i$.
\end{proof}

\begin{lem} \label{genfinitecor}
Let $f: X \to Z$ be a surjective generically finite map from a smooth projective variety $X$ to a $\mathbb{Q}$-factorial normal projective variety $Z$.  Let $L$ be an $\mathbb{R}$-Cartier divisor on $X$ and let $\{ T_{i} \}_{i = 1}^{k}$ be a collection of irreducible curves on $X$.  Suppose that there are constants $\alpha_{i}$ with
\begin{equation*}
L \cdot C = (\deg f|_{C}) \alpha_{i}
\end{equation*}
for every curve $C$ on $X$ with $f(C) = f(T_{i})$.
\begin{enumerate}
\item Suppose that for each $i$ the image $f(T_{i})$ is a curve lying in the open locus on $Z$ over which $f$ is flat.  Then there is an $\mathbb{R}$-Cartier divisor $D$ on $Z$ such that $L \cdot T_{i} = f^{*}D \cdot T_{i}$ for every $i$.
\item Suppose that the numerical classes of the $T_{i}$ span $N_{1}(X)$ and that $L \cdot C = 0$ for every $f$-vertical curve $C$.  Then $L \equiv f^{*}D$ for some $\mathbb{R}$-Cartier divisor $D$ on $Z$.
\end{enumerate}
\end{lem}

\begin{proof}
Choose a normal birational model $\phi: X' \to X$ and a normal birational model $\mu: Z' \to Z$ so that we have a morphism $f': X' \to Z'$ flattening $f$.  We may assume that $\phi$ and $\mu$ are isomorphisms on the locus over which $f$ is flat.  For each $i$ choose an irreducible curve $T_{i}'$ on $X'$ lying above $T_{i}$. 

We first prove (1).  Suppose that $C$ is a curve on $X'$ such that $f'(C) = f'(T_{i}')$.  Since $f(\phi(T_{i}'))$ is a curve,
\begin{align*}
\phi^{*}L \cdot C & = (\deg \phi|_{C}) (\deg f|_{\phi(C)}) \alpha_{i} \\
& = (\deg f'|_{C})(\deg \mu|_{f'(C)}) \alpha_{i} \\
& = (\deg f'|_{C})  (\deg \mu|_{f'(T_{i}')})\alpha_{i}.
\end{align*}
Set $\alpha_{i}' = (\deg \mu|_{f'(T_{i}')}) \alpha_{i}$.  The set of curves $\{ T_{i}' \}$ and the divisor $\phi^{*}L$ satisfy the hypotheses of Lemma \ref{finitepullback} for the finite map $f'$ and the constants $\alpha_{i}'$.  Lemma \ref{finitepullback} yields a divisor $D_{Z'}$ on $Z'$ such that $\phi^{*}L \cdot T_{i}' = f'^{*}D_{Z'} \cdot T_{i}'$ for every $i$.  Since the $T_{i}$ lie over the locus on which $f$ is flat, $f'(T_{i}')$ avoids the $\mu$-exceptional locus.  Thus, setting $D = \mu_{*}D_{Z'}$, we obtain $L \cdot T_{i} = f^{*}D \cdot T_{i}$ for every $i$.

We next prove (2).  Applying Lemma \ref{birationallem} to $\phi$, we can find finitely many irreducible $\phi$-vertical curves $\{ S_{j}' \}_{j =1}^{r}$so that the span of the numerical classes of the $T_{i}'$ and $S_{j}'$ is all of $N_{1}(X')$.

Suppose that $C$ is a curve on $X'$ such that $f'(C) = f'(T_{i}')$.  If $f(\phi(T_{i}'))$ is a curve, then as before set $\alpha_{j}' =  (\deg \mu|_{f'(T_{j}')})\alpha_{j}$.  If $f(\phi(T_{i}'))$ is a point, then $C$ is also $(\mu \circ f')$-vertical.  Since $\phi^{*}L$ has vanishing intersection with every $\mu \circ f'$-vertical curve, $\phi^{*}L \cdot C = 0 = \phi^{*}L \cdot T_{i}'$.  Set $\alpha_{i}' = 0$.  Similarly, set $\beta_{j}' = 0$ for every $S_{j}'$.

The set of curves $\{ T_{i}' \} \cup \{ S_{j}' \}$ and the divisor $\phi^{*}L$ satisfy the hypotheses of  Lemma \ref{finitepullback} for the finite map $f'$ and the constants $\alpha_{i}',\beta_{j}'$.  The result of the lemma indicates that there is a divisor $D_{Z'}$ on $Z'$ so that $f'^{*}D_{Z'} \equiv \phi^{*}L$.  Since $D_{Z'}$ is $\mu$-numerically trivial, Lemma \ref{birationallem} yields a divisor $D$ on $Z$ such that $\mu^{*}D \equiv D_{Z'}$.  Thus $f^{*}D \equiv L$.
\end{proof}

\section{Numerical Triviality on Every Fiber} \label{everyfibersection}

In this section we give the proof of Theorem \ref{numtrivialeveryfiber}.  We start by recalling an example demonstrating that the pseudo-effectiveness of $L$ is necessary in order to have any hope of relating $L$ to divisors on the base.

\begin{exmple}
Let $E$ be an elliptic curve without complex multiplication and consider the surface $S = E \times E$ with first projection $\pi: S \to E$.  Recall that $N^{1}(S)$ is generated by the fibers $F_{1},F_{2}$ of the two projections and the diagonal $\Delta$.  In particular, the subspace of $\pi$-trivial divisors is generated by $F_{1}$ and $\Delta - F_{2}$.  Since this space has larger dimension than $N^{1}(E)$, most $\pi$-trivial divisors will not be numerically equivalent to pull-backs from $E$.
\end{exmple}

The first step in the proof is to show that numerical equivalence of divisors can be detected against curves which are intersections of very ample divisors.

\begin{prop} \label{amplespan}
Suppose that $X$ is a normal projective variety of dimension $n$.  Fix a set of ample Cartier divisors $\mathcal{H} = \{ H_{1},\ldots,H_{r} \}$ whose numerical classes span $N^{1}(X)$.  Then the intersection products $\mathcal{H}^{n-1}$ span $N_{1}(X)$.
\end{prop}

\begin{proof}
We first show that any irreducible curve $C$ on $X$ is contained in an irreducible surface $S$ that is the complete intersection of members of the linear systems $|m_{i}H_{i}|$ for some sufficiently large $m_{i}$.  For a fixed curve $C$ and ample divisor $H \in \mathcal{H}$, choose $m$ sufficiently large so that $h^{1}(X,\mathcal{I}_{C} \otimes \mathcal{O}_{X}(mH)) = 0$ and $\mathcal{O}_{X}(mH)$ is globally generated.  Then $\mathcal{I}_{C} \otimes \mathcal{O}_{X}(mH)$ is globally generated away from $C$.  Furthermore, for sufficiently large $m$ the corresponding linear system is not composite with a pencil.  Thus, the Bertini theorems show that the general member of this linear system is irreducible and is normal away from $C$.  Arguing inductively, we construct the surface $S$.

Now suppose that $L$ is a divisor on $X$ such that $L \cdot \alpha = 0$ for every element $\alpha$ in the span of $\mathcal{H}^{n-1}$.  Consider an irreducible curve $C$ on $X$ and the corresponding irreducible surface $S$.  We know that $L|_{S} \cdot H|_{S} = 0$ for any ample divisor $H \in \mathcal{H}$.  Furthermore $L|_{S}^{2} = 0$ since the $H_{i}$ span $N^{1}(X)$.  The Hodge Index Theorem (applied to a resolution of $S$) implies that $L|_{S} \equiv 0$, and in particular, $L \cdot C =0$.  Since $C$ is arbitrary, $L$ is numerically trivial.
\end{proof}

We now turn to the proof of Theorem \ref{numtrivialeveryfiber}.

\begin{proof}[Proof of \ref{numtrivialeveryfiber}:]
Suppose that $\phi: X' \to X$ is a resolution of $X$.  Note that $\phi^{*}L$ is $(f \circ \phi)$-numerically trivial.  If we can find a $D$ so that $\phi^{*}L \equiv (f \circ \phi)^{*}D$ then also $L \equiv f^{*}D$.  Thus we may assume that $X$ is smooth.

We next show that for a curve $R$ through a very general point of $Z$ there is some constant $\alpha_{R}$ such that
\begin{equation*} \label{earlycomp} \tag{$\dagger$}
L \cdot C = (\deg f|_{C}) \alpha_{R}
\end{equation*}
for every curve $C$ on $X$ with $f(C) = R$.  Let $R'$ denote the normalization of $R$ and consider the normalization $Y$ of $X \times_{Z} R'$.  Since $R$ goes through a very general point of $Z$ we may assume that the pullback of $L$ to every component of $Y$ is pseudo-effective and that only one component of $Y$ dominates $S'$.  Consider two curves $C$ and $C'$ on $Y$ with $f(C) = f(C') = R'$.  By cutting down $Y$ by very ample divisors, we can find a chain of normal surfaces $S_{i}$ connecting $C$ to $C'$, all of which map surjectively to $R'$ under $f$.  We may ensure that $L|_{S_{i}}$ is pseudo-effective for every $i$.

Applying Corollary \ref{numtrivsurfaces} to the surface $S_{i}$, we see that there is some divisor $D_{i}$ on $R'$ such that $L|_{S_{i}} \equiv f^{*}D_{i}$.  For $i \geq 1$ let $C'_{i}$ denote the curve $S_{i} \cap S_{i+1}$.  Since $C'_{i}$ dominates $R'$, we have
\begin{equation*}
\deg(D_{i}) \deg(f|_{C'_{i}}) = L \cdot C'_{i} = \deg(D_{i+1}) \deg(f|_{C'_{i}}).
\end{equation*}
Thus there is one fixed $D_{1}$ so that $L|_{S_{i}} \equiv f^{*}D_{1}$ for every $i$.  Fixing $C$ and letting $C'$ vary, we see that the constant $\alpha_{R} = \deg(D_{1})$ satisfies the desired condition for every curve above $R$.

Let $W \subset X$ denote a smooth very general intersection of very ample divisors such that the map $f: W \to Z$ is generically finite and $L|_{W}$ is pseudo-effective.  Certainly $L|_{W}$ has vanishing intersection with any $f$-vertical curve on $W$.  Furthermore, Proposition \ref{amplespan} shows that we can choose a finite collection of curves $T_{i}$ through very general points whose numerical classes span $N_{1}(X)$.  In particular \eqref{earlycomp} holds over the $T_{i}$.  Lemma \ref{genfinitecor} (2) yields a divisor $D$ on $Z$ such that $L|_{W} \equiv f^{*}D$. 

Apply Proposition \ref{amplespan} to $X$ to find a collection of irreducible curves $C_{i}$ on $X$ that are not $f$-vertical and whose numerical classes span $N_{1}(X)$.  For each $i$ choose an irreducible curve $C_{i}^{W}$ on $W$ such that $f(C_{i}) = f(C_{i}^{W})$.  Since
\begin{align*}
L \cdot C_{i} & = (\deg f|_{C_{i}})\alpha_{f(C_{i})} \\
& = (\deg f|_{C_{i}}) \frac{L \cdot C_{i}^{W}}{\deg f|_{C_{i}^{W}}} \\
& = (f^{*}D \cdot C_{i}^{W}) \frac{\deg f|_{C_{i}}}{\deg f|_{C_{i}^{W}}} \\
& = f^{*}D \cdot C_{i}
\end{align*}
we see that $L \equiv f^{*}D$.
\end{proof}

\section{Numerical Triviality on a General Fiber} \label{generalfibersection}

In this section we prove Theorem \ref{numtrivialgeneralfiber}.  In fact, we will address the more general situation where $\nu(L|_{F})=0$ for a general fiber $F$.  The following examples show that Theorem \ref{numtrivialgeneralfiber} is optimal in some sense.

\begin{exmple}
Let $f:S \to C$ be a morphism from a smooth surface to a smooth curve.  Suppose that $L$ is an effective $f$-degenerate divisor.  Then $L$ is not numerically equivalent to the pull-back of a divisor on the base.  This is still true on higher birational models of $f$.  One must pass to the positive part $P_{\sigma}(L) = 0$.
\end{exmple}

\begin{exmple}
Let $D$ be a big divisor on a smooth variety $X$ and let $\phi: Y \to X$ be a blow-up along a smooth center along which $D$ has positive asymptotic valuation.  Then $P_{\sigma}(\phi^{*}D) < \phi^{*}P_{\sigma}(D)$ is not numerically equivalent to a pull-back of a divisor on $X$.  One must pass to the flattening $id: Y \to Y$.
\end{exmple}

The following is a stronger version of Theorem \ref{numtrivialgeneralfiber}.

\begin{thrm}
Let $f: X \to Z$ be a surjective morphism of normal projective varieties with connected fibers.  Suppose that $L$ is a pseudo-effective $\mathbb{R}$-Cartier divisor such that $\nu(L|_{F}) = 0$ for a general fiber $F$ of $f$.  Then there is a smooth birational model $\phi: Y \to X$, a map $g: Y \to Z'$ birationally equivalent to $f$, and an $\mathbb{R}$-Cartier divisor $D$ on $Z'$ such that $P_{\sigma}(\phi^{*}L) \equiv P_{\sigma}(g^{*}D)$.
\end{thrm}

\begin{proof}
By passing to a resolution we may assume that $X$ is smooth.

Choose a normal birational model $\mu: X' \to X$, a smooth birational model $Z'$ of $Z$, and a morphism $f': X' \to Z'$ flattening $f$.  Let $\psi: Y \to X'$ denote a smooth model.  We let $g$ denote the composition $f' \circ \psi$ and let $\phi$ denote the composition $\mu \circ \psi$.  Note that every $g$-exceptional divisor is also $\phi$-exceptional.

Since $\nu(\phi^{*}L|_{F})=0$ for a general fiber $F$ of $g$, we have
\begin{equation*}
\phi^{*}L|_{F} \equiv N_{\sigma}(\phi^{*}L|_{F}) \leq N_{\sigma}(\phi^{*}L)|_{F}.
\end{equation*}
In particular $P_{\sigma}(\phi^{*}L)|_{F} \equiv 0$ for a general fiber $F$. 

By Corollary \ref{psigmarelativelytrivial}, there is a subset $V \subset Z'$ that is a countable union of codimension $2$ subsets such that $P_{\sigma}(\phi^{*}L)$ is numerically trivial along every fiber not over $V$.  In particular, suppose that the curve $R \subset Z'$ avoids $V$ and $P_{\sigma}(\phi^{*}L)$ is pseudo-effective when restricted to the fiber over $R$.  By the same argument as in the proof of Theorem \ref{numtrivialeveryfiber}, there is some constant $\alpha_{R}$ such that
\begin{equation*} \label{compare} \tag{*}
P_{\sigma}(\phi^{*}L) \cdot C = \deg(g|_{C}) \cdot \alpha_{R}
\end{equation*}
for every curve $C$ with $g(C)=R$.

We next apply the generically finite case to construct a divisor $D_{1}$.  Choose a smooth very general intersection $W$ of very ample divisors on the smooth variety $Y$ so that the map $g|_{W}: W \to Z$ is generically finite.  By choosing $W$ very generically we may assume that the divisor $P_{\sigma}(\phi^{*}L)|_{W}$ is pseudo-effective.

Consider the subspace of $N_{1}(W)$ generated by irreducible curves $C$ that avoid $g^{-1}(V)$ and run through a very general point of $W$.  We may choose a finite collection of irreducible curves $\{ T_{i} \}$ satisyfing these two properties whose numerical classes span this subspace. Thus there are constants $\alpha_{i}$ so that
\begin{equation*} 
P_{\sigma}(\phi^{*}L)|_{W} \cdot C = \deg(g|_{C}) \cdot \alpha_{i}
\end{equation*}
for every curve $C$ with $g(C)=g(T_{i})$.   Applying Lemma \ref{genfinitecor} (1), we find a divisor $D_{1}$ on $Z'$ with $P_{\sigma}(\phi^{*}L) \cdot T_{i} = g^{*}D_{1} \cdot T_{i}$ for every $i$.  Furthermore $P_{\sigma}(\phi^{*}L) \cdot C = g^{*}D_{1} \cdot C$ for any curve $C$ through a very general point of $W$ such that $g(C)$ avoids $V$, since $C$ is numerically equivalent to a sum of the $T_{i}$.

We next relate $P_{\sigma}(\phi^{*}L)$ and $g^{*}D_{1}$.  Recall that $\mu(f'^{-1}V)$ is a countable union of codimension $2$ subvarieties in $X$.  By Proposition \ref{amplespan} we may choose curves $S_{i}^{X}$ avoiding this locus and running through a very general point of $X$ whose numerical classes form a basis for $N_{1}(X)$.  Let $\{ S_{i} \}$ consist of the strict transforms of these curves on $Y$.  Since the $S_{i}$ are generic, for each we may choose a curve $S_{i}^{W} \subset W$ going through a very general point and such that $g(S_{i}^{W}) = g(S_{i})$ avoids $V$.  This guarantees that $L \cdot S_{i}^{W} = g^{*}D_{1} \cdot S_{i}^{W}$.  By construction $P_{\sigma}(\phi^{*}L) \cdot S_{i}^{W}$ and $P_{\sigma}(\phi^{*}L) \cdot S_{i}$ can be compared using \eqref{compare}.  Arguing as in the proof of Theorem \ref{numtrivialeveryfiber}, we see that $P_{\sigma}(\phi^{*}L) \cdot S_{i} = g^{*}D_{1} \cdot S_{i}$ for every $i$.  This proves that
\begin{equation*}
\phi_{*}P_{\sigma}(\phi^{*}L) \equiv \phi_{*}g^{*}D_{1}.
\end{equation*}
Choose effective $\phi$-exceptional divisors $E$ and $F$ with no common components such that
\begin{equation*}
P_{\sigma}(\phi^{*}L) + E \equiv g^{*}D_{1} + F.
\end{equation*}
Note that since $f: X \to Z$ is generically flat, no $\phi$-exceptional divisor dominates $Z'$.  In particular $F$ is $g$-vertical, so we may apply Lemma \ref{makefdegen} to $F$ to find
\begin{equation*}
F = g^{*}D_{2} + F_{deg} - F_{exc}
\end{equation*}
where $F_{deg}$ is $g$-degenerate and $F_{exc}$ is $g$-exceptional.  Set $D = D_{1} + D_{2}$. Then
\begin{equation*} \label{psig} \tag{**}
P_{\sigma}(\phi^{*}L) + E + F_{exc} \equiv g^{*}D + F_{\deg}.
\end{equation*}
Since $F_{deg}$ is $g$-degenerate, Lemma \ref{fdegen} shows $P_{\sigma}(g^{*}D + F_{deg}) = P_{\sigma}(g^{*}D)$.  Similarly, since $(E+F_{exc})$ is $\phi$-exceptional,
\begin{align*}
P_{\sigma}(\phi^{*}L) & \leq P_{\sigma}(P_{\sigma}(\phi^{*}L)+E+F_{exc}) \\
& \leq P_{\sigma}(\phi^{*}L + E + F_{exc}) \\
& = P_{\sigma}(\phi^{*}L) \textrm{ by Lemma \ref{fdegen}.}
\end{align*}
Taking the positive part of both sides of  \eqref{psig} yields $P_{\sigma}(\phi^{*}L) \equiv P_{\sigma}(g^{*}D)$.
\end{proof}

\nocite{*}
\bibliographystyle{amsalpha}
\bibliography{trivmorph}

\newcommand{\etalchar}[1]{$^{#1}$}
\providecommand{\bysame}{\leavevmode\hbox to3em{\hrulefill}\thinspace}
\providecommand{\MR}{\relax\ifhmode\unskip\space\fi MR }
\providecommand{\MRhref}[2]{%
  \href{http://www.ams.org/mathscinet-getitem?mr=#1}{#2}
}
\providecommand{\href}[2]{#2}
\begin{thebibliography}{BCE{\etalchar{+}}02}

\bibitem[BCE{\etalchar{+}}02]{8authors}
T.~Bauer, F.~Campana, T.~Eckl, S.~Kebekus, T.~Peternell, S.~Rams, T.~Szemberg,
  and L.~Woltzlas, \emph{A reduction map for nef line bundles}, Complex
  Geometry (G{\"o}ttingen, 2000), Springer, Berlin, 2002, pp.~27--36.

\bibitem[BDPP04]{bdpp04}
S.~Boucksom, J.P. Demailly, M.~P{\v a}un, and T.~Peternell, \emph{The
  pseudo-effective cone of a compact {K\"a}hler manifold and varieties of
  negative {K}odaira dimension}, 2004,
  http://www-fourier.ujf-grenoble.fr/{$\sim$}demailly/manuscripts/coneduality.pdf,
  submitted to J. Alg. Geometry.

\bibitem[Bou04]{boucksom04}
S.~Boucksom, \emph{Divisorial {Z}ariski decompositions on compact complex
  manifolds}, Ann. Sci. {\'E}cole Norm. Sup. \textbf{37} (2004), no.~4, 45--76.

\bibitem[Eck05]{eckl05}
T.~Eckl, \emph{Numerically trival foliations, {I}itaka fibrations, and the
  numerical dimension}, 2005, arXiv:math/0508340v1.

\bibitem[Fuj79]{fujita79}
T.~Fujita, \emph{On {Z}ariski problem}, Proc. Japan Acad. Ser. A Math. Sci.
  \textbf{55} (1979), no.~3, 106--110.

\bibitem[GL11]{gl11}
Y.~Gongyo and B.~Lehmann, \emph{Reduction maps and minimal model theory}, 2011,
  arXiv:1103.1605v1.

\bibitem[Leh11]{lehmann11}
B.~Lehmann, \emph{On {E}ckl's pseudo-effective reduction map}, 2011,
  arXiv:1103.1073v1 [math.AG].

\bibitem[Nak04]{nakayama04}
N.~Nakayama, \emph{Zariski-decomposition and abundance}, MSJ Memoirs, vol.~14,
  Mathematical Society of Japan, Tokyo, 2004.

\bibitem[Pet08]{peternell08}
T.~Peternell, \emph{Varieties with generically nef tangent bundles}, 2008,
  arXiv:0807.0982v1 [math.AG].

\bibitem[Ray72]{raynaud72}
M.~Raynaud, \emph{Flat modules in algebraic geometry}, Comp. Math. \textbf{24}
  (1972), 11--31.

\bibitem[Zar64]{zariski64}
O.~Zariski, \emph{The theorem of {R}iemann-{R}och for high multiples of an
  effective divisor on an algebraic surface}, Ann. of Math. \textbf{76} (1964),
  no.~2, 560--615.

\end{thebibliography}

\end{document}